\documentclass{article}
\usepackage{amsmath}
\usepackage{amsfonts}
\usepackage{amsthm}
\theoremstyle{plain}
\newtheorem{theorem}{Theorem}[section]
\newtheorem{lemma}[theorem]{Lemma}

\theoremstyle{definition}

\theoremstyle{remark}

\begin{document}
\title{Convergence in capacity\\
}
\author{Urban Cegrell}

\maketitle

\section{Introduction}
The purpose of this paper is to study convergence of Monge-Amp\`ere
measures associated to sequences of plurisubharmonic functions defined
on a hyperconvex subdomain $\Omega$ of $\mathbb C^n$. 

This paper is an updated version of the preprint $[C0]$.

The concept of convergence in capacity was introduced in $[X]$, where
it was proved that for any uniformly bounded sequence $\varphi_{j}$ of
plurisubharmonic functions that converges in $\mathbb C^n$-capacity we
have that $(dd^c\varphi_{j})^n$ converges weak$^{*}$
to
$(dd^c\varphi)^n, j\to+\infty$.

We generalize this result:

\begin{theorem}
    Assume $\mathcal F\ni u_{0}\le u_{j}\in\mathcal F$ and that $u_{j}$
converges to $u$
    in $\mathbb C^n$-capacity. Then $(dd^cu_{j})^n$ converges weak* to
    $(dd^cu)^n, j\to+\infty$.
    \end{theorem}

    We first recall some definitions. See
    [C1] and [C2] for details.

    The class $\mathcal F$ consists of all plurisubharmonic functions
$\varphi$
    on $\Omega$ such that there is a sequence $\varphi_{j}\in\mathcal E_{0},
    \varphi_{j}\searrow\varphi, j\to+\infty$ and $\sup\limits_{j}
    \int\limits_{\Omega}(dd^c\varphi_{j})^n<+\infty$, where $\mathcal
E_{0}$ is
    the class of bounded plurisubharmonic functions $\psi$ such that
    $\lim\limits_{z\to\xi}\psi(z)=0, \forall\xi\in\partial\Omega$ and
    $\int\limits_{\Omega}(dd^c\psi)^n<+\infty$.

    The following definition was introduced in [X]: A sequence
    $\varphi_{j}\in\mathcal F$ converges to $\varphi$ in $\mathbb
    C^n$-capacity if
    $$\mathrm{cap}(\{z\in k:|\varphi-\varphi_{j}|>\delta\})\to0,\
    j\to+\infty\quad \forall\ k\subset\subset\Omega,\ \forall\
    \delta>0.$$
    For $\omega\subset\subset\Omega,
    \mathrm{cap}(\omega)=\int\limits_{\Omega}(dd^ch_{\omega}^{*})^n$ where
    $h^{*}_{\omega}$ is the smallest upper semicontinuous majorant of
    $$h_{\omega}(z)=
    \sup\{\varphi(z); \varphi\in \mathcal E_0; \varphi|_{\omega}\le -1\}.$$
    Finally, we write $u_{j}\rightsquigarrow u,\ j\to+\infty$ if $u_{j}$
    converges weak* to $u, j\to+\infty$.

    \newpage

    \section{Proofs}

    \begin{lemma}
Suppose $\mu$ is a positive measure on $\Omega$ which vanishes on
all pluripolar
sets and $\mu(\Omega)<+\infty$. If $u_{0}, u_{j}\in\mathcal F,\
u_{0}\le
u_{j}\rightsquigarrow u,\ j\to+\infty$ and if $\int
u_{0}d\mu>-\infty$, then $\lim\limits_{j\to+\infty}\int
u_{j}d\mu=\int ud\mu$.
    \end{lemma}

     \begin{proof} Denote by $dV$ the Lebesgue measure and choose
     $\tilde u_{j}\in\mathcal E_{0}\cap C(\tilde \Omega),\ \tilde u_{j}\ge
     u_{j}$ such that $\int\limits_{\Omega}(\tilde
     u_{j}-u_{j})(d\mu+dV)<\frac{1}{j}$. Then \hbox{$\tilde
     u\rightsquigarrow u,\ \lim\limits_{j\to+\infty}\int
     u_{j}d\mu-\int\tilde u_{j}d\mu=0$} so it is enough to prove that
     $$\lim\limits_{j\to+\infty}\int\tilde u_{j}d\mu=\int ud\mu.$$
     Thus we can assume $u_{j}\in\mathcal E_{0}\cap C(\bar \Omega)$.

     By Theorem 6.3 in [C1] there is a $\psi\in\mathcal E_{0}, \ f\in
     L^1\big((dd^c\psi)^n\big)$ with
     $$\mu=f(dd^c\psi)^n$$
     so by lemma 5:2 in [C1], for every $p<+\infty$,
     $$\lim\limits_{j\to+\infty}\int u_{j}d\mu_{p}=\int ud\mu_{p}$$
     there $\mu_{p}=\min(f,p)(dd^c\psi)^n$.

     Now
     \begin{equation*}
     \begin{split}
     &\lim\limits_{j\to+\infty}\int
u_{j}d\mu=\lim\limits_{j\to+\infty}\int u_{j}d\mu_{p}+\\
&+\lim\limits_{j\to+\infty}\int
u_{j}\big(f-\mathrm{min}(f,p)\big)(dd^c\psi)^n\ge\\
&\ge\int ud\mu_{p}+\int
u_{0}\big(f-\mathrm{min}(f,p)\big)d\mu\\
&\rightarrow\int ud\mu,\ p\to+\infty
\end{split}
\end{equation*}
by monotone convergence. On the other hand, by Fatous lemma,
$$\limsup\limits_{j\to+\infty}\int u_{j}d\mu\le \int ud\mu$$
which gives the desired conclusion.
\end{proof}

\begin{proof}
Proof of the theorem.

We prove that
$$\lim\limits_{j\to+\infty}\int h(dd^cu_{j})^n=
\int h(dd^cu)^n,\
\forall\ h\in\mathcal E_{0},$$
which is enough by Lemma 3:1 in [C1].

Suppose
$\omega_{1},\omega_{2},\dots,\omega_{n-1}\in\mathcal F,\
h\in\mathcal E_{0}$.
It follows
from the assumption that $u_{j}\rightsquigarrow u,\
j\to+\infty$ so
$\lim\limits_{j}\int \omega_{1}*\ dd^cu_{j}\land
dd^c\omega_{2}\land\dots\land
dd^ch=\int\omega_{1}dd^cu\land\dots\land dd^ch$ by the lemma.

Suppose now that
$$\lim\limits_{j}\int\omega_{1}(dd^cu_{j})^p\land
dd^c\omega_{p+1}\land\dots\land dd^ch=
\int\omega_{1}(dd^cu)^p\land\dots\land dd^ch.$$
for $1\leq q\leq p\leq n-2.$
We claim
\begin{equation*}
\begin{split}
&\lim\limits_{j}\int\omega_{1}(dd^cu_{j})^{p+1}\land dd^c\omega_{p+2}
\land\dots\land dd^c\omega_{n-1}\land dd^ch=\\
&=\int\omega_{1}(dd^cu)^{p+1}\land
dd^c\omega_{p+2}\land\dots
\land dd^c\omega_{n-1}\land dd^ch.
\end{split}
\end{equation*}
Given $\varepsilon>0$ choose $k\subset\subset\Omega$ such that
$\{z\in\Omega ; h<-\varepsilon\}\subset k$ and then a subsequence
$u_{j_{t}}$ such that
$$\int\limits_{\Omega}\bigg(dd^c\sum\limits^\infty_{t=1}h^{*}_{\{z\in
k ; |u-u_{j_{t}}|>\varepsilon\}}\bigg)^n<1$$
and denote by
$$h_{N}=\max\bigg(\sum\limits^\infty_{t=N}h^{*}_{\{z\in
k ; |u-u_{j_{t}}|>\varepsilon\}}, -1\bigg).$$
Then $h_{N}\to 0, N\to+\infty$ outside a pluripolar set and
\begin{equation*}
    \begin{split}
    &\int\omega_{1}(dd^cu_{j_{k}})^{p+1}\land\dots\land
    dd^ch-\int\omega_{1}(dd^cu_{j_{k}})^{p}dd^cu\land\dots\land dd^ch=\\
    &=\int(u_{j_{k}}-u)\big(dd^cu_{j_{k}}\big)^p\land
    dd^c\omega_{1}\land\dots\land dd^ch=\\
    &=\int-h_{N}(u_{j_{k}}-u)\big(dd^cu_{j_{k}}\big)^p\land
    dd^c\omega_{1}\land\dots\land dd^ch+\\
    &+\int(1+h_{N})(u_{j_{k}}-u)\big(dd^cu_{j_{k}}\big)^p\land
    dd^c\omega_{1}\land\dots\land dd^c(h-h_{\varepsilon})+\\
    &+\int(1+h_{N})(u_{j_{k}}-u)\big(dd^cu_{j_{k}}\big)^p\land
    dd^c\omega_{1}\land\dots\land dd^ch_{\varepsilon}
    =I_{j_{k}}+II_{j_{k}}+III_{j_{k}}.
    \end{split}
    \end{equation*}

where
$h_{\varepsilon}=\max(h,-\varepsilon).$

Thus
\begin{equation*}
    \begin{split}
    &|I_{j_{k}}|\le
    2\int(-h_{N})(-u_{0})(dd^cu_{j_{k}})^p\land\dots\land
    dd^ch\le\\
&\le\int\big[-u_{0}+\max(u_{0_{1}}-R)\big](dd^cu_{j_{k}})^p\land\dots\land
dd^ch+\\
    &+R\int-h_{N}\big(dd^cu_{j_{k}}\big)^p\land\dots\land dd^ch\\
    &\to\int\big[-u_{0}+\max(u_{0},-R)\big](dd^cu)^p\land\dots\land dd^ch+\\
    &+R\int(-h_{N})(dd^cu)^p\land\dots\land dd^ch,\ j_{k}\to+\infty.
    \end{split}
    \end{equation*}
Now, the first term on the right hand side is small if $R$ is large,
the second is small if $N$ is even larger.

\begin{equation*}
    \begin{split}
&|III_{j_{k}}|\le \int -u_{0}*(dd^cu_{j_{k}})^p\land
dd^c\omega_{1}\land\dots\land dd^ch_{\varepsilon}\le\\
&\le \int -u_{0}*(dd^cu_{0})^{p}\land dd^c\omega_{1}\land\dots\land
dd^ch_{\varepsilon} \le\\
&\le \int -h_{\varepsilon}(dd^cu_{0})^{p+1}\land dd^c\omega_{1}
\land\dots\land dd^c\omega_{n-1}\le\\
&\varepsilon\int(dd^cu_{0})^{p+1}\land dd^c\omega_{1}\dots\land
dd^c\omega_{n-1}*
\end{split}
\end{equation*}
so it remains to estimate
\begin{equation*}
    \begin{split}
&|II_{j_{k}}|\le 2\varepsilon \bigg(\int (dd^cu_{j_{k}})^p\land
dd^c\omega_{1}\dots \land dd^c(h+h_{\varepsilon})\bigg)\le\\
&\le 4\varepsilon\int (dd^cu_{0})^p\land dd^c\omega_{1}\dots\land
dd^ch.
\end{split}
\end{equation*}

    which proves the claim.

    Using the claim and repeating the chain of inequalities above, we
    conclude that
    \begin{equation*}
    \begin{split}
     &\int h(dd^cu_{j_{k}*}*)^n - \int h(dd^cu)^n =  \\
&\int (u_{j_{k}*}* - u)(dd^cu_{j_{k}*}*)^{n-1}\land h  + \\
&\int u(dd^cu_{j_{k}*}*)^{n-1}\land h - \int u(dd^cu)^{n-1}\land h \to 0,
j_{k}*\ \to +\infty.
\end{split}
\end{equation*}

    This proves the theorem, since we have proved that every
    subsequence of
    $u_{j}*$ contains a subsequence $u_{j_{t}*}*$ such that
    $(dd^cu_{j_{t}*}*)^n$ converges weak* to $(dd^cu)^n, t\to+\infty .$

    \end{proof}

\today


\begin{thebibliography}{9}
\bibitem[C0]] Cegrell, U., {\it Convergence in capacity.}
Isaac Newton Institute for Mathematical Sciences. Preprint Series
NI01046-NPD. Cambridge. 2001.

\bibitem[C1]] Cegrell, U., {\it Pluricomplex energy.} Acta Math.
{\bf 180:2} (1998), 187--217.

\bibitem[C2]] Cegrell, U.,  {\it The general definition of the
complex Monge-Amp\`ere operator.} Ann. Inst. Fourier 54 (2004).

\bibitem[X]] Xing, Y., {\it Continuity of the complex
Monge-Amp\`ere operator.}
Proc. AMS., {\bf 124:2} (1996), 457--467.

\end{thebibliography}
 \end{document}